\documentclass[12pt]{article}
\usepackage{graphicx} 
\usepackage{amsmath}
\usepackage{amssymb}
\begin{document}
\newcommand{\vare}{\varepsilon} 
\newtheorem{tth}{Theorem}[section]
\newtheorem{dfn}[tth]{Definition}
\newtheorem{lem}[tth]{Lemma}
\newtheorem{prop}[tth]{Proposition}
\newtheorem{coro}[tth]{Corollary}
\renewcommand{\quad}{\hspace*{1em}}
\begin{center}
{\Large {\bf Genericity of Caustics on a corner} }
\vspace*{0.4cm}\\
{\large Takaharu Tsukada}\footnote{Higashijujo 3-1-16 
Kita-ku, Tokyo 114-0001
JAPAN. e-mail : tsukada@math.chs.nihon-u.ac.jp}
\vspace*{0.2cm}\\
{\large  College of Humanities \& Sciences, Department of Mathematics,\\
 Nihon University}
\end{center}
\begin{abstract}
We introduce the notions of  {\em the caustic-equivalence} and {\em the weak caustic-equivalence relations} of reticular Lagrangian maps
 in order to give a generic classification of  caustics on a corner. 
We give the figures of all generic caustics on a corner in a smooth manifold of 
dimension $2$ and $3$.
\end{abstract}

\section{Introduction}
\quad
In \cite{retLag} 
we investigate the theory of {\em reticular Lagrangian maps} which can be 
described stable caustics generated by a hypersurface germ with an $r$-corner in a smooth manifold.
A map germ  
\[ \pi \circ i:({\mathbb L},0) \rightarrow (T^* {\mathbb R}^n,0) \rightarrow({\mathbb R}^n,0)\]
is called {\em a reticular Lagrangian map} if $i$ is a restriction of a symplectic diffeomorphism germ
on $(T^* {\mathbb R}^n,0)$,
where $I_r=\{1,\ldots,r\}$ and 
${\mathbb L}=\{(q,p)\in T^* {\mathbb R}^n |q_1p_1=\cdots 
=q_rp_r=q_{r+1}=\cdots=q_n=0,q_{I_r}\geq 0  \}$
be a representative of the union of
\[L_\sigma^0=
\{(q,p)\in (T^* {\mathbb R}^n,0)|q_\sigma=p_{I_r-\sigma}=q_{r+1}=\cdots=q_n=0, 
q_{I_r-\sigma}\geq 0 \}  
\mbox{ for all }\sigma\subset I_r.\]
We define the caustic of $\pi\circ i$ is the union of 
the caustics $C_\sigma$ of the Lagrangian maps $\pi\circ i |_{L_\sigma^0} $for all $\sigma\subset I_r$ 
and the quasi-caustic $Q_{\sigma,\tau}=\pi\circ i(L_\sigma^0\cap L_\tau^0)$
 for all $\sigma,\tau\subset I_r(\sigma \neq \tau)$.
In the case $r=2$, that is the initial hypersurface germ has a corner, the caustic of $\pi\circ i$ is 
\[ C_\emptyset\cup C_1\cup C_2 \cup C_{1,2} \cup Q_{\emptyset,1} \cup Q_{\emptyset,2} 
\cup Q_{1,\{1,2\}} \cup Q_{2,\{1,2\}}. \]
For the definitions of generating families of reticular Lagrangian maps,
see \cite[p.575-577]{retLag}.
In \cite{generic}
we investigate the genericity of caustics on an $r$-corner
and give the generic classification for the cases $r=0$ and $1$ by using 
 G.Ishikawa's methods (see \cite[Section 5]{ishikawa3}).
We also showed that this method do not work well for the case
$r=2$.
In this paper we introduce the two equivalence relations of reticular Lagrangian maps
which are weaker than Lagrangian equivalence in order to
give a generic classification of caustics on a corner.\\
\begin{figure}[ht]
\begin{center}
Example of stable caustics on a corner \\
 \begin{minipage}{0.40\hsize}
  \begin{center}
    \includegraphics[width=5.5cm,height=5cm]{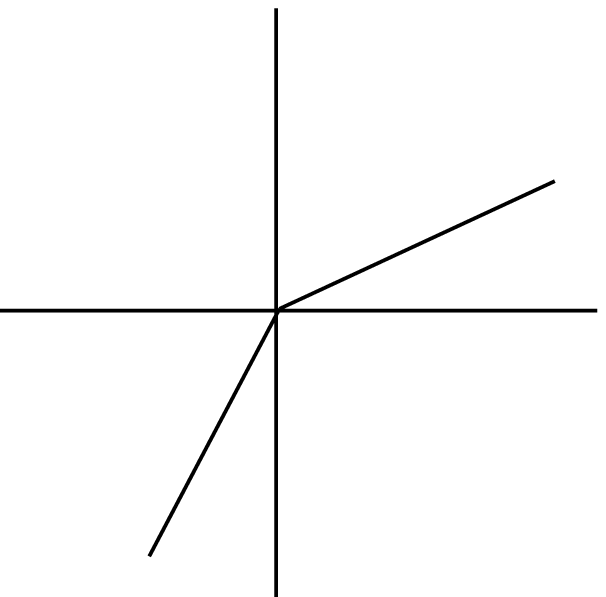}
  \caption{weakly caustic-stable}  
\end{center}
 \end{minipage}\hspace{2cm}
 \begin{minipage}{0.40\hsize}
  \begin{center}
    \includegraphics[width=5.5cm,height=5cm]{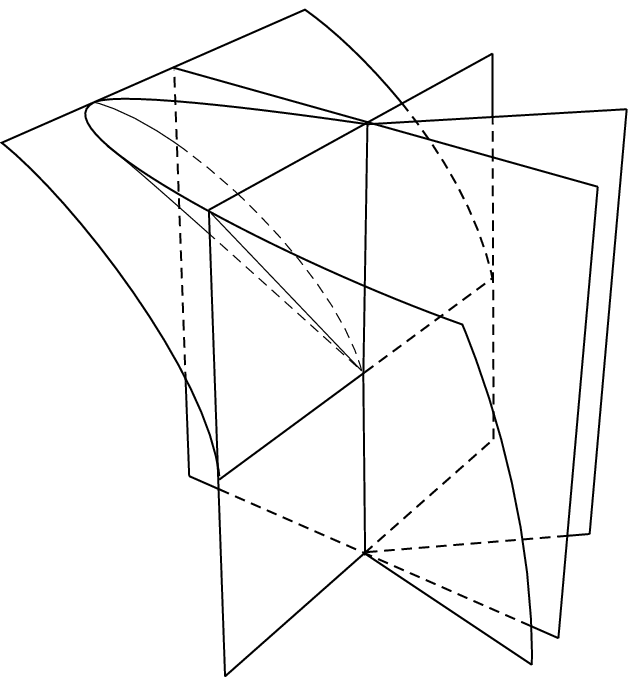}
\caption{caustic-stable}
\end{center}
 \end{minipage}
\end{center}
\end{figure}

\section{Caustic-equivalence and Weak caustic-equivalence}\label{caust:sec}
\quad
We introduce the equivalence relations of reticular Lagrangian maps and 
their generating families.

Let $\pi\circ i_j$ be reticular Lagrangian maps for $j=1,2$.
We say that  they are {\em caustic-equivalent} if there exists a diffeomorphism 
germ $g$ on $({\mathbb R}^n,0)$ such that 
\begin{equation}
 g(C^1_\sigma)=C^2_\sigma,\ \ g(Q^1_{\sigma,\tau})=Q^2_{\sigma,\tau}\ \ 
\mbox{ for all } \sigma,\tau \subset I_r\ (\sigma \neq \tau). \label{caust:def}
\end{equation}

We say that  reticular Lagrangian maps $\pi\circ i_1$ and $\pi\circ i_2$ are {\em weakly caustic-equivalent} 
if there exists a homeomorphism 
germ $g$ on $({\mathbb R}^n,0)$ such that 
 $g$ is smooth on all $C^1_\sigma$, $Q^1_{\sigma,\tau}$, and satisfies (\ref{caust:def}).\\

We shall define the stabilities of reticular Lagrangian maps under the above equivalence relations
and define the corresponding equivalence relations and stabilities of their generating families. \\

The purpose of this paper is to show the following theorem:
\begin{tth}\label{genericclassLag:th}
Let $n=2,3$ or $4$, $U$ a neighborhood of $0$ in $T^* {\mathbb R}^n$,
$S(T^*{\mathbb R}^n,0)$ be the set of symplectic diffeomorphism germs on $(T^*{\mathbb R}^n,0)$,
and $S(U,T^*{\mathbb R}^n)$ be the space of symplectic embeddings from 
$U$ to $T^*{\mathbb R}^n$ with $C^\infty$-topology.
Then there exists a residual set $O\subset S(U,T^* {\mathbb R}^n)$
such that for any $\tilde{S}\in O$ and $x\in U$,
the reticular Lagrangian map $\pi\circ \tilde{S}_x|_{\mathbb L}$ 
is weakly caustic-stable or caustic-stable, 
where $\tilde{S}_x\in S(T^* {\mathbb R}^n,0)$ be defined by the map $x_0\mapsto \tilde{S}(x_0+x)-\tilde{S}(x)$.
\end{tth}
A reticular Lagrangian map $\pi\circ \tilde{S}_x|_{\mathbb{L}}$ for any $\tilde{S}\in O$ and $x\in U$ is 
weakly caustic-equivalent to one which has a generating family
$B_{2,2}^{\pm,+,1},B_{2,2}^{\pm,+,2},
B_{2,2}^{\pm,-}$, 
or is caustic equivalent to one which has a generating family 
$B_{2,2}^{\pm,0},B_{2,2,3}^{\pm.\pm},B_{2,3}^{\pm,\pm},
B_{3,2}^{\pm,\pm},B_{2,3'}^{\pm,\pm},B_{3,2'}^{\pm,\pm},C_{2,3}^{\pm,\pm},C_{3,2,1}^{\pm,\pm}
,C_{3,2,2}^{\pm,\pm}$.\\
$B_{2,2}^{\pm,+,1}$: $F(x_1,x_2,q_1,q_2)=x_1^2\pm x_1x_2+\frac15x_2^2+q_1x_1+q_2x_2$,\\
$B_{2,2}^{\pm,+,2}$: $F(x_1,x_2,q_1,q_2)=x_1^2\pm x_1x_2+x_2^2+q_1x_1+q_2x_2$,\\
$B_{2,2}^{\pm,-}$: $F(x_1,x_2,q_1,q_2)=x_1^2\pm x_1x_2-x_2^2+q_1x_1+q_2x_2$,\\
$B_{2,2}^{\pm,0}$: $F(x_1,x_2,q_1,q_2,q_3)=x_1^2\pm x_2^2+q_1x_1+q_2x_2+q_3x_1x_2$,\\
$B_{2,2,3}^{\pm,\pm}$: $F(x_1,x_2,q_1,q_2,q_3)=(x_1\pm x_2)^2\pm x_2^3+q_1x_1+q_2x_2+q_3x_2^2$,\\
$B_{2,3}^{\pm,\pm}$: $F(x_1,x_2,q_1,q_2,q_3)=x_1^2\pm x_1x_2\pm x_2^3+q_1x_1+q_2x_2+q_3x_2^2$,\\
$B_{3,2}^{\pm,\pm}$: $F(x_1,x_2,q_1,q_2,q_3)=x_1^3\pm x_1x_2\pm x_2^2+q_1x_1+q_2x_2+q_3x_1^2$,\\
$B_{2,3'}^{\pm,\pm}$: $F(x_1,x_2,q_1,q_2,q_3,q_4)=x_1^2\pm x_1x_2^2\pm  x_2^3
+q_1x_2^2+q_2x_1x_2+q_3x_2+q_4x_1$,\\
$B_{3,2'}^{\pm,\pm}$: $F(x_1,x_2,q_1,q_2,q_3,q_4)=x_1^3\pm x_1^2x_2\pm  x_2^2
+q_1x_1^2+q_2x_1x_2+q_3x_1+q_4x_2$,\\
$C_{3,2}^{\pm,\pm}$: $F(y,x_1,x_2,q_1,q_2,q_3)=\pm y^3+x_1y\pm x_2y +x_2^2+q_1y+q_2x_1+q_3x_2$,\\
$C_{3,2,1}^{\pm,\pm}$: $F(y,x_1,x_2,q_1,q_2,q_3,q_4)=\pm y^3+ x_1y\pm  x_2y^2+x_2^2
+q_1y^2+q_2y+q_3x_1x_2+q_4x_2$,\\
$C_{3,2,2}^{\pm,\pm}$: $F(y,x_1,x_2,q_1,q_2,q_3,q_4)=\pm y^3+ x_2y\pm  x_1y^2+x_1^2
+q_1y^2+q_2y+q_3x_1x_2+q_4x_1$.\\

In order to describe the caustic-equivalence of reticular Lagrangian maps by
their generating families, we introduce the following equivalence relation of function germs.
We say that function germs $f,g\in{\cal E}(r;k)$ are {\em reticular ${\cal C}$-equivalent} if
there exist $\phi\in{\cal B}(r;k)$ and a non-zero number $a\in {\mathbb R}$ 
such that $g=a \cdot f\circ \phi$. 
See \cite{retLag} or \cite{tPKfunct} for the notations.
We construct the theory of unfoldings with respect to the corresponding equivalence relation.
Then the relation of unfoldings is given as follows:
Two function germs $F(x,y,q),G(x,y,q)\in {\cal E}(r;k+n)$ are
{\em reticular ${\cal P}$-${\cal C}$-equivalent} if
there exist $\Phi\in{\cal B}_n(r;k+n)$ 
and a unit $a\in {\cal E}(n)$ and $b \in {\cal E}(n)$ and  such that $G=a\cdot F\circ \Phi+b$. 
We define the {\em stable} reticular (${\cal P}$-)${\cal C}$-equivalence by the ordinary ways
(see \cite[p.576]{retLag}).
We remark that a reticular ${\cal P}$-${\cal C}$-equivalence class includes the reticular 
${\cal P}$-${\cal  R}^+$-equivalence classes.\\

%
%
We review the results of the theory.
Let $F(x,y,u)\in {\mathfrak M}(r;k+n)$ be an unfolding of 
$f(x,y)\in {\mathfrak M}(r;k)$.

%
We say that $F$ is {\em reticular ${\cal P}$-${\cal C}$-stable} if the following condition holds: For any neighborhood $U$ of
$0$ in ${\mathbb R}^{r+k+n}$ and any representative $\tilde{F} \in C^\infty 
(U,{\mathbb R})$ of $F$, 
there exists a neighborhood $N_{\tilde{F}}$ of $\tilde{F}$ in $C^\infty$-topology
 such that for any element $\tilde{G} \in N_{\tilde{F}}$
the germ $\tilde{G}|_{{\mathbb H}^r\times {\mathbb R}^{k+n}}$ at 
$(0,y_0,q_0)$ is reticular ${\cal P}$-${\cal C}$-equivalent to $F$
for some $(0,y_0,q_0)\in U$.\\

%
We say that $F$ is {\em reticular ${\cal P}$-${\cal C}$-versal} if all unfolding of $f$ is  reticular ${\cal P}$-${\cal C}$-$f$-induced from $F$.
That is, for any unfolding $G\in {\mathfrak M}(r;k+n')$ of $f$, there exist 
$\Phi\in{\mathfrak M}(r;k+n',r;k+n)$ and
a unit $a \in {\cal E}(n')$ and $b\in {\cal E}(n')$
satisfying the following conditions: \\
(1) $\Phi(x,y,0)=(x,y,0)$ for all $(x,y)\in ({\mathbb H}^r\times {\mathbb R}^k,0)$ and $a(0)=1,\ b(0)=0$,\\
(2) $\Phi$ can be written in the form:
\[ \Phi(x,y,q)=(x_1\phi^1_1(x,y,q),\cdots,x_r\phi^r_1(x,y,q),
\phi_2(x,y,q),\phi_3(q)),\]
(3) $G(x,y,q)=a(q)\cdot F\circ\Phi(x,y,q)+b(q)$ for all 
    $(x,y,q)\in ({\mathbb H}^r\times {\mathbb R}^{k+n'},0)$.\\

%
We say that $F$ is {\em reticular ${\cal P}$-${\cal C}$-infinitesimally versal} if 
\[ 
{\cal E}(r;k)=\langle x \frac{\partial f}{\partial x},
\frac{\partial f}{\partial
y}\rangle_{ {\cal E}(r;k) }+
\langle 1,f,\frac{\partial F}{\partial q}|_{q=0}\rangle_{{\mathbb R}}.
\]

%
We say that $F$ is {\em reticular ${\cal P}$-${\cal C}$-infinitesimally stable} if
\[ {\cal E}(r;k+n)  = \langle x \frac{\partial
 F}{\partial x},\frac{\partial
 F}{\partial y}\rangle_{ {\cal
 E}(r;k+n) }+\langle 1,F,\frac{\partial F}{\partial q}\rangle_{{\cal E}(n)}. 
\]
\vspace{3mm}

%
\quad 
We say that $F$ is {\em reticular ${\cal P}$-${\cal C}$-homotopically stable} if for any
smooth path-germ $({\mathbb R},0)\rightarrow {\cal E}(r;k+n),t\mapsto F_t$ with
$F_0=F$, there exists a smooth path-germ $({\mathbb R},0)\rightarrow {\cal
B}_n(r;k+n)\times {\cal E}(n)\times {\cal E}(n),t\mapsto (\Phi_t,a_t,b_t)$ with
$(\Phi_0,a_0,b_0)=(id,1,0)$ such that each $(\Phi_t,a_t,b_t)$ is a reticular ${\cal P}$-${\cal C}$-isomorphism from $F$ to $F_t$, that is
$F_t=a_t\cdot F\circ \Phi_t+b_t$ for $t$ around $0$.

%
%
\begin{tth}{\rm (cf., \cite[Theorem 4.5]{retLag})}
Let $F\in {\mathfrak M}(r;k+n)$ be an unfolding of $f\in {\mathfrak M}(r;k)$.
Then the following are all equivalent. \\
{\rm (1)} $F$ is reticular ${\cal P}$-${\cal C}$-stable.\\
{\rm (2)} $F$ is reticular ${\cal P}$-${\cal C}$-versal.\\
{\rm (3)} $F$ is reticular ${\cal P}$-${\cal C}$-infinitesimally versal. \\
{\rm (4)} $F$ is reticular ${\cal P}$-${\cal C}$-infinitesimally stable. \\
{\rm (5)} $F$ is reticular ${\cal P}$-${\cal C}$-homotopically stable. 
\end{tth}

For a non-quasihomogeneous function germ $f(x,y) \in {\mathfrak M}(r;k)$, 
if $1,f,a_1,\ldots,a_n\in {\cal E}(r;k)$ is a representative of a basis of the 
vector space 
\[ {\cal E}(r;k)/ \langle x \frac{\partial f}{\partial x},
\frac{\partial f}{\partial y}\rangle_{ {\cal E}(r;k)}, \]
then 
the function germ $f+a_1q_1+\cdots +a_nq_n\in {\mathfrak M}(r;k+n)$ is a 
reticular ${\cal P}$-${\cal C}$-stable unfolding of $f$.
We call $n$ the reticular ${\cal C}$-codimension of $f$.
If $f$ is quasihomogeneous then 
$f$ is included in $\langle x \frac{\partial f}{\partial x},
\frac{\partial f}{\partial y}\rangle_{ {\cal E}(r;k)}$.
This means that the reticular ${\cal C}$-codimension of a quasihomogeneous function germ 
 is equal to its reticular ${\cal R}^+$-codimension.

We define the {\em simplicity} of function germs under the reticular ${\cal C}$-equivalence
in the usual way (cf., \cite{generic}). 
\begin{tth}\label{simpleC:th}{\rm (cf., \cite[Theorem 2.1,2.2]{generic})}
 A reticular ${\cal C}$-simple function germ in ${\mathfrak M}(1;k)^2$ is 
stably reticular ${\cal C}$-equivalent to one of the following function germs:
\[
B_l:x^l\  (l\geq 2),\ \ \ C^{\vare}_l:xy+\vare y^{l}\ (\vare^{l-1}=1,
l\geq 3),\ \ \ 
F_4: x^2+y^3.\]
 \end{tth} 
 
The relation between reticular Lagrangian maps and 
their generating families under the caustic-equivalence are given as follows:
\begin{prop}\label{wcequi:prop}
Let $\pi\circ i_j$ be reticular Lagrangian maps with generating families $F_j$ 
for $j=1,2$.
If $F_1$ and $F_2$ are stably reticular ${\cal P}$-${\cal C}$-equivalent then 
$\pi\circ i_1$ and $\pi\circ i_2$ are caustic-equivalent.
\end{prop}
{\em Proof.} The function germ $F_2$ may be written that $F_2(x,y,q)=a(q)F_3(x,y,q)$,
where $a$ is a unit and $F_1$ and $F_3$ are stably reticular ${\cal P}$-${\cal R}^+$-equivalent.
Then the reticular Lagrangian map $\pi\circ i_3$ given by $F_3$ and $\pi\circ i_1$ are 
Lagrangian equivalent and the caustic of $\pi\circ i_2$ and $\pi\circ i_3$ coincide to each other.
\hfill $\blacksquare$

This proposition shows that it is enough to classify function germs under 
the stable  reticular ${\cal P}$-${\cal C}$-equivalence in order to classify reticular Lagrangian maps
under the caustic-equivalence.
We here give the classification list as the following:
\begin{tth}{\rm (cf., \cite[p.592]{retLag})}
Let $f\in {\mathfrak M}(2;k)^2$ have the reticular ${\cal C}$-codimension$\leq 4$.
Then $f$ is stably reticular ${\cal C}$-equivalent to one of the following list.
\begin{center}
\begin{tabular}{lllll}
\hline $k$ & Normal form & codim & Conditions & Notation \\ 
\hline
 $0$ & $x_1^2\pm x_1x_2+a x_2^2$ & $3$ & $0<a<\frac{1}{4}$ & 
                              $B_{2,2,a}^{\pm,+,1}$\\
 & $x_1^2\pm x_1x_2+a x_2^2$ & $3$ & $a>\frac{1}{4}$ & 
                              $B_{2,2,a}^{\pm,+,2}$\\
 & $x_1^2\pm x_1x_2+a x_2^2$ & $3$ & $a<0$ & 
                              $B_{2,2,a}^{\pm,-}$\\
 & $x_1^2\pm x_2^2$ & $3$ &  & 
                              $B_{2,2}^{\pm,0}$\\
 & $(x_1\pm x_2)^2\pm x_2^3$ & $3$ & & 
                             $B_{2,2,3}^{\pm,\pm}$\\
& $x_1^2\pm x_1x_2\pm x_2^3$ & $3$  & &
                              $B_{2,3}^{\pm,\pm}$\\
 & $x_1^3\pm x_1x_2\pm  x_2^2$ & $3$ & & 
                              $B_{3,2}^{\pm,\pm}$ \\
 & $x_1^2\pm x_1x_2^2\pm  x_2^3$ & $4$ & & 
                              $B_{2,3'}^{\pm,\pm}$ \\
 & $x_1^3\pm x_1^2x_2\pm  x_2^2$ & $4$ & & 
                              $B_{3,2'}^{\pm,\pm}$ \\                          
\hline 
 $1$ & $\pm y^3+ x_1y\pm  x_2y+x_2^2$ & $3$ 
      
      &  & $C_{3,2}^{\pm,\pm}$ \\
  & $\pm y^3+ x_1y\pm  x_2y^2+x_2^2$ & $4$ 
      
      &  & $C_{3,2,1}^{\pm,\pm}$ \\
  & $\pm y^3+ x_2y\pm  x_1y^2+x_1^2$ & $4$ 
      
      &  & $C_{3,2,2}^{\pm,\pm}$ \\
      
\hline 
\end{tabular}
\end{center}
\end{tth}
We remark that the stable reticular ${\cal C}$-equivalence class $B_{2,3}^{+,+}$ of 
$x_1^2+x_1x_2+x_2^3$ consists of the union of the stable reticular ${\cal R}$-equivalence classes
of $x_1^2+x_1x_2+ax_2^3$ and $-x_1^2-x_1x_2-ax_2^3$ for $a>0$.
The same things hold for $B_{2,2,3}^{\pm,\pm}$, $B_{2,3}^{\pm,\pm}$, $B_{3,2}^{\pm,\pm}$, $C_{3,2}^{\pm,\pm}$.

%
%
%
%
%
\section{Caustic-stability}
\quad
We define {\em the caustic-stability} of reticular Lagrangian maps and reduce our investigation 
to finite dimensional jet spaces of symplectic diffeomorphism germs.

We say that a reticular Lagrangian map $\pi\circ i$ is {\em caustic-stable} if 
the following condition holds:
For any extension $S\in S(T^*{\mathbb R}^n,0)$ of $i$ and any representative 
$\tilde{S}\in S(U,T^*{\mathbb R}^n)$ of $S$, there exists a neighborhood $N_{{\tilde{S}}}$ of 
$\tilde{S}$ such that for any $\tilde{S}'\in N_{{\tilde{S}}}$ the reticular Lagrangian map
$\pi\circ \tilde{S}'|_{\mathbb L}$ at $x_0$ and $\pi\circ i$ are caustic-equivalent for some 
$x_0=(0,\ldots,0,p^0_{r+1},\ldots,p^0_n)\in U$.
\begin{dfn}\label{caustequi:dfn}
Let $\pi\circ i$ be a reticular Lagrangian map
and $l$ be a non-negative number.
We say that $\pi\circ i$ is {\em caustic $l$-determined} if the following condition holds:
For any extension $S$ of $i$, the reticular Lagrangian map 
$\pi\circ S'|_{{\mathbb L}}$ and $\pi\circ i$ are caustic-equivalent for any symplectic diffeomorphism germ $S'$ on 
$(T^* {\mathbb R}^n,0)$ satisfying $j^lS(0)=j^lS'(0)$.
\end{dfn}
\begin{lem}\label{findet:lem}
Let $\pi\circ i:({\mathbb L},0) \rightarrow (T^* {\mathbb R}^n,0) 
\rightarrow({\mathbb R}^n,0)$ be a reticular Lagrangian map.
If a generating family of $\pi\circ i$ is 
reticular ${\cal P}$-${\cal C}$-stable then $\pi\circ i$ is caustic $(n+2)$-determined.
\end{lem}
{\em Proof}. This is proved by the analogous method of \cite[Theorem 5.3]{generic}.
We give the sketch of proof.
Let $S$ be an extension of $i$.
Then we may assume that there exists a function germ $H(Q,p)$ such that the
canonical relation $P_S$ has the form:
\[ P_S=\{(Q,-\frac{\partial H}{\partial Q}(Q,p),
-\frac{\partial H}{\partial p}(Q,p),p)\in 
(T^* {\mathbb R}^n\times T^* {\mathbb R}^n,(0,0))
\}.\]
Then the function germ $F(x,y,q)=H_0(x,y)+\langle y,q\rangle$
is a reticular ${\cal P}$-${\cal C}$-stable generating family of $\pi\circ i$,
and $H_0$ is reticular ${\cal R}$-$(n+3)$-determined, 
where $H_0(x,y)=H(x,0,y)$.
Let a symplectic diffeomorphism germ $S'$ on $(T^* {\mathbb R}^n,0)$ 
satisfying $j^{n+2}S(0)=j^{n+2}S'(0)$ be given.
Then there exists a function germ $H'(Q,p)$ such that 
the canonical relation $P_{S'}$ is given the same form for $H'$ and
the function germ $G(x,y,q)=H'_0(x,y)+\langle y,q\rangle$ is a 
generating family of $\pi\circ S'|_{{\mathbb L}}$.
Then it holds that $j^{n+3}H_0(0)=j^{n+3}H_0'(0)$.
There exists a function germ $G'$ such that $G$ and $G'$ are
 reticular ${\cal P}$-${\cal R}$-equivalent and   
$F$ and $G'$ are reticular ${\cal P}$-${\cal C}$-infinitesimal versal 
unfoldings of $H_0(x,y)$.
It follows that $F$ and $G$ are reticular ${\cal P}$-${\cal C}$-equivalent.
Therefore $\pi\circ i$ and $\pi \circ S'|_{{\mathbb L}}$ are
caustic-equivalent.\hfill $\blacksquare$\\

For a reticular ${\cal P}$-${\cal C}$-stable unfolding $F\in {\mathfrak M}(2;k+n)^2$ with $n\leq 3$, 
the function germ $f=F|_{q=0}$ has a modality under the reticular ${\cal R}$-equivalence
(see \cite[p.592]{retLag}).
For example, consider the case $f$ is stably reticular ${\cal C}$-equivalent to 
$x_1^2+x_1x_2+x_2^3$.
Then $F$ is stably reticular ${\cal P}$-${\cal C}$-equivalent to $f+q_1x_1+q_2x_2+q_3x_2^2$.
In this case the function germs $F_a(x,q)=x_1^2+x_1x_2+ax_2^3+q_1x_1+q_2x_2+q_3x_2^2( a>0)$
are stably reticular ${\cal P}$-${\cal C}$-equivalent to $F$ but not 
stably reticular ${\cal P}$-${\cal R}^+$-equivalent to each other.
Let $S_a^\pm $ be extensions of reticular Lagrangian embeddings defined by $F_a$ and
 $-F_a$ for $a>0$ respectively.
We define the caustic-equivalence class of $S_1$ by 
\[ [S_1]_c:=\bigcup_{a>0}([S^+_a]_L\cup [S^-_a]_L),\]
where $[S^\pm_a]_L$ are the Lagrangian equivalence classes of $S^\pm_a$ respectively.
By Proposition \ref{wcequi:prop}, we have that 
all reticular Lagrangian maps $\pi\circ S'|_{{\mathbb L}}$ are caustic-equivalent to 
each other for $S'\in [S_1]_c$.\\

In order to apply the transversality theorem to our theory, we need to prove that 
the set consists of the $2$-jets of the caustic-equivalence class $[S_1]_c$,
we denote this by $[j^2S_1(0)]_c$, is an immersed manifold of $S^2(3)$,
where $S^l(n)$ be the smooth manifold which consists of $l$-jets of elements in $S(T^*{\mathbb R}^n,0)$.
We shall prove that the map germ $(0,\infty)\rightarrow S^2(3),a\mapsto j^2S_a(0)$
is not tangent to $[j^2S_a(0)]_L$ for any $a$, and apply the following lemma:
\begin{lem}\label{moduli:tth}
Let $I$ be an open interval, $N$ a manifold, and $G$ a Lie group acts on $N$ smoothly.
Suppose the orbits $G\cdot x$ have the same dimension for all $x\in I$.
Let $x:I\rightarrow N$ be a smooth path such that 
$\frac{dx}{dt}(t)$ is not tangent to $G\cdot x(t)$
for all $t\in I$. Then
\[ \bigcup_{t\in I} G\cdot x(t) \]
is an immersed manifold of $N$.
\end{lem}
We denote that we here prove the case $B_{2,3}^{+,+}$. The same method is valid for all $B_{2,3}^{\pm,\pm},
B_{3,2}^{\pm,\pm}$.

We define $G_a\in {\mathfrak M}(6)^2$ by $G_a(Q_1,Q_2,Q_3,q_1,q_2,q_3)=F_a(Q_1,Q_2,q_1,q_2)+Q_3q_3$.
Then $G_a$ define the canonical relations $P_a$ and they give symplectic diffeomorphisms $S_a$ of the forms:
\[ S_a(Q,P)=(-2Q_1-Q_2-P_1,-Q_1-3aQ_2^2-P_2+2P_3Q_2,-P_3,Q_1,Q_2,Q_2^2+Q_3).\] 
We have that $F_a$ are generating families of $\pi\circ S_a|_{{\mathbb L}}$.
Then $\frac{d S_a}{d a}=(0,-3Q_2^2,0,0,0,0)=X_f\circ S_a$ for $f=-p_2^3$.
We suppose that $j^2(\frac{d S_a}{d a})(0)\in T_z([z]_L)$ for $z=j^2S_a(0)$.
By \cite[Lemma 6.2]{generic}, 
there exist a fiber preserving function germ $H\in {\mathfrak M}_{Q,P}^2$ and 
$g\in \langle Q_1P_1,Q_2P_2 \rangle_{{\cal E}_{Q,P}}+{{\mathfrak M}_{Q,P}}\langle Q_3\rangle$ 
such that $j^2(X_f\circ S_a)(0)=j^2(X_H\circ S_a+(S_a)_*X_g)(0)$.
This means that $j^3(f\circ S_a)(0)=j^3(H\circ S_a+g)(0)$.
It follows that 
there exist function germs $h_1,h_2,h_3\in {\mathfrak M}_Q$, $h_0\in {\mathfrak M}^2_Q$ such that
\begin{eqnarray*}
 f\circ S_a=-Q_2^3 & \equiv & h_1(q\circ S_a)Q_1+h_2(q\circ S_a)Q_2+h_3(q\circ S_a)(Q_2^2+Q_3)+h_0(q\circ S_a) \\
 & & \mbox{  mod }
\langle Q_1P_1,Q_2P_2 \rangle_{{\cal E}_{Q,P}} +{{\mathfrak M}_{Q,P}}\langle Q_3 \rangle +{\mathfrak M}^4_{Q,P}.
\end{eqnarray*}
We may reduce this to
\begin{eqnarray*}
-Q_2^3 & \equiv & h_1(-2Q_1-Q_2,-Q_1-3aQ_2^2-P_2+2P_3Q_2,-P_3)Q_1\\
& & +h_2(-2Q_1-Q_2-P_1,-Q_1-3aQ_2^2+2P_3Q_2,-P_3)Q_2\\
& & +h_3(-2Q_1-Q_2-P_1,-Q_1,-P_3)Q_2^2+h_0(-2Q_1-Q_2-P_1,-Q_1-P_2,-P_3)\\
& &  \mbox{  mod }
\langle Q_1P_1,Q_2P_2 \rangle_{{\cal E}_{Q,P}} +{\mathfrak M}^{}_{Q,P}\langle Q_3 \rangle +{\mathfrak M}^4_{Q,P}.
\end{eqnarray*}
We show this equation has a contradiction.
The coefficients of $P_1^{i_1}P_2^{i_2}P_3^{i_3}$ on the equation depend only on the coefficients of 
$q_1^{i_1}q_2^{i_2}q_3^{i_3}$ on $h_0$
respectively.
This means that $h_0(q\circ S_a) \equiv 0$.
The coefficients of $Q_1^2,Q_1P_2,Q_1P_3$ on the equation depend only on 
the coefficients of $q_1,q_2,q_3$ on $h_1$ respectively.
This means that $j^1(h_1(q\circ S_a)(0)\equiv 0$.
The coefficients of $Q_2P_1,Q_1Q_2,Q_2P_3$ on the equation depend only on 
the coefficients of $q_1,q_1,q_3$ on $h_2$.
This means that $j^1(h_2(q\circ S_a))(0) \equiv 0$.
So we need only to consider the quadratic part of $h_1,h_2$ and the linear part of $h_3$.
The coefficients of $Q_2P_1^2,Q_2^2P_1$ on the equation depend only on
the coefficient of $q_1^2$ on $h_2$ and the coefficient of $q_1$ on $h_3$ respectively. 
This means that their coefficients are all equal to $0$.
Therefore the coefficient of $Q_2^3$ on the right hand side of the equation is $0$.
This contradicts the equation. 
So we have that $j^2(\frac{d S_a}{d a})(0)$ is not included in  $T_z([z]_L)$.\\

We also prove the case $B_{2,2,3}^{+,+}$: We consider the reticular Lagrangian maps $\pi\circ i_a$ with the generating families 
$F_a(x_1,x_2,q_1,q_2,q_3)=(x_1+x_2)^2+ax_2^3+q_1x_1+q_2x_2+q_3x_2^2$.
Then the function germs
$G_a(Q_1,Q_2,Q_3,q_1,q_2,q_3)=(Q_1+Q_2)^2+aQ_2^3+q_1Q_1+q_2Q_2+q_3Q_2^2+q_3Q_3$ 
are the generating functions of the canonical relations $P_{S_a}$ and
$i_a=S_a|_{\mathbb{L}}$.
Then $S_a$ have the forms:
\[ S_a(Q,P)=(-(2Q_1+2Q_2+P_1),-(2Q_1+2Q_2+3aQ_2^2+P_2-2P_3Q_2),-P_3,Q_1,Q_2,Q_2^2+Q_3).\]
We have that  $\frac{d S_a}{d a}=(0,-3Q_2^2,0,0,0,0)=X_f\circ S_a$ for $f=-p_2^3$.
Then we consider the following equation:
\begin{eqnarray*}
 f\circ S_a=-Q_2^3 & \equiv & h_1(q\circ S_a)Q_1+h_2(q\circ S_a)Q_2+h_3(q\circ S_a)(Q_2^2+Q_3)+h_0(q\circ S_a) \\ 
& & \mbox{  mod }
\langle Q_1P_1,Q_2P_2 \rangle_{{\cal E}_{Q,P}} +{{\mathfrak M}_{Q,P}}\langle Q_3 \rangle +{\mathfrak M}^4_{Q,P},
\end{eqnarray*}
where $h_1,h_2,h_3\in {\mathfrak M}(Q), h_0\in {\mathfrak M}^2(Q)$.
We may reduce this to
\begin{eqnarray*}
-Q_2^3& \equiv & h_1(-(2Q_1+2Q_2),-(2Q_1+2Q_2+3Q_2^2+P_2-2Q_2P_3),-P_3)Q_1\\
& & +
h_2(-(2Q_1+2Q_2+P_1),-(2Q_1+2Q_2+3Q_2^2-2Q_2P_3),-P_3)Q_2\\
& & +h_3(-(2Q_1+2Q_2+P_1),-(2Q_1+2Q_2),-P_3)Q_2^2\\
& & + h_0(-(2Q_1+2Q_2+P_1),-(2Q_1+2Q_2+3aQ_2^2+P_2-2Q_2P_3),-P_3) \\
& & \mbox{  mod }
\langle Q_1P_1,Q_2P_2 \rangle_{{\cal E}_{Q,P}} +{{\mathfrak M}_{Q,P}}\langle Q_3 \rangle +{\mathfrak M}^4_{Q,P}.
\end{eqnarray*}
By the same reason in the case $B_{2,3}^{+,+}$, we have that $h_0(q\circ S_a)\equiv 0$.
By the consideration of the coefficients of $Q_1^2,Q_1P_2,Q_1P_3$  and
$Q_2P_1,Q_2^2,Q_2P_3$ on the equation, we have that $j^1(h_1(q\circ S_a)Q_1)(0)\equiv
j^1(h_2(q\circ S_a)Q_2)(0)\equiv 0$.
The coefficients of $Q_1P_2^2,Q_1P_3^2,Q_1P_2P_3$ on the equation depend only on
the coefficients of $q_2^2,q_3^2,q_2q_3$ on $h_1$.
This means that they are all equal to $0$.
The coefficients of $Q_1^2P_2,Q_1^2P_3,Q_1^3$ depend only on the coefficients of $q_1q_2,q_1q_3,q_1^2$
on $h_1$.
This means that they are all equal to $0$. We have that 
$j^2(h_1(q\circ S_a)Q_1)(0)\equiv 0$.

The coefficients of $Q_2P_1^2,Q_2P_3^2,Q_2P_1P_3$ depend only on the coefficients of $q_1^2,q_3^2,q_1q_3$
on $h_2$ and  they are all equal to $0$.
We  write $h_2=q_2(bq_1+cq_2+dq_3), h_3=eq_1+fq_2+gq_3$.
We calculate the coefficients of $Q_1^2Q_2,Q_1Q_2^2,Q_2^2P_1,Q_1Q_2P_3,Q_2^2P_3$, 
then we have that $-2 b - 2 c=-8 (-2 b - 2 c) + 2 e (-2 - 2 f)=4 b - 2 e=d=4 d - 2 e g=0$.
This is solved that $b=c=d=e=0$ or $b=\frac{e}2,c=-\frac{e}2,d=0,f=-1,g=0$.
This means that the coefficient of $Q_2^3$ on the right hand side of the equation is 
$4 b + 4 c - 2 e - 2 e f=0$.
This contradicts the equation.\\

We also prove the case $C_{3,2}^{+,+}$: 
We consider the reticular Lagrangian maps $\pi\circ i_a$ with the generating families 
$F_a(y,x_1,x_2,q_1,q_2,q_3)=y^3+x_1y+x_2y+ax_2^2+ax_2^3+q_1y+q_2x_1+q_3x_2$.
Then the function germs
$G_a(y,Q_1,Q_2,Q_3,q_1,q_2,q_3)=y^3+Q_1y+Q_2y+aQ_2^2+q_1y+q_2Q_1+q_3Q_2+yQ_3$ 
are the generating families of the canonical relations $P_{S_a}$ and
$i_a=S_a|_{\mathbb{L}}$.
Then $S_a$ have the forms:
\[ S_a(Q,P)=(-(3P_3^2+Q_1+Q_2+Q_3),P_3-P_1,P_3-2aQ_2-P_2,-P_3,Q_1,Q_2).\]
We have that  $\frac{d S_a}{d a}=(0,0,-2Q_2,0,0,0)=X_f\circ S_a$ for $f=-p_3^2$.
Then we consider the following equation:
\begin{eqnarray*}
 f\circ S_a=-Q_2^2 & \equiv & h_1(q\circ S_a)(-P_3)+h_2(q\circ S_a)Q_1+h_3(q\circ S_a)Q_2+h_0(q\circ S_a) \\
 & & \mbox{  mod }
\langle Q_1P_1,Q_2P_2 \rangle_{{\cal E}_{Q,P}} +{{\mathfrak M}_{Q,P}}\langle Q_3 \rangle +{\mathfrak M}^3_{Q,P}.
\end{eqnarray*}
We may reduce this to
\begin{eqnarray*}
-Q_2^2 & \equiv & h_1(-(Q_1+Q_2),P_3-P_1,P_3-2aQ_2-P_2)(-P_3)\\
& & +h_2(-(Q_1+Q_2),P_3,P_3-2aQ_2-P_2)Q_1\\
& & +h_3(-(Q_1+Q_2),P_3-P_1,P_3-2aQ_2)Q_2\\
& & +h_0(-(Q_1+Q_2),P_3-P_1,P_3-2aQ_2-P_2)\\
& &  \mbox{  mod }
\langle Q_1P_1,Q_2P_2 \rangle_{{\cal E}_{Q,P}} +{{\mathfrak M}_{Q,P}} \langle Q_3 \rangle+{\mathfrak M}^3_{Q,P}.
\end{eqnarray*}
Since the coefficients of $P_1^{i_2}P_2^{i_3}$ on the equation depend only on 
the coefficients of $q_2^{i_2}q_3^{i_3}$ on $h_0$, it follows that they are all equal to $0$.
Since the coefficients of $P_1P_3,P_2P_3$ depend only on 
the coefficients of $q_2,q_3$ on $h_1$, it follows that they are all equal to $0$.

Therefore we may set $h_1=bq_1,\ h_2=cq_1+dq_2+eq_3,\ 
h_3=fq_1+gq_2+hq_3, \ 
h_0=q_1(iq_1+jq_2+hq_3)$.
By the calculation of the equation, we have that
 the coefficient of $Q_2^2$ on the right hand side of the equation is $0$.
This contradicts the equation.

\begin{lem}
Let $\pi\circ i:({\mathbb L},0) \rightarrow (T^* {\mathbb R}^n,0) 
\rightarrow({\mathbb R}^n,0)$ be a reticular Lagrangian map,
$S$ be an extension of $i$.
Suppose that 
the caustic-equivalence class $[j^{n+2}_0S(0)]_c$ be an immersed manifold of $S^{n+2}(n)$.
If a generating family of $\pi\circ i$ is reticular ${\cal P}$-${\cal C}$-stable and
$j^{n+2}_0S$ is transversal to $[j^{n+2}_0S(0)]_c$ at $0$,
then $\pi\circ i$ is caustic stable.
\end{lem}
This is proved by the analogous method of \cite[Theorem 6.6 (t)\&(is)$\Rightarrow$(s)]{generic}.
By this lemma, we have that the caustic-stability of reticular  Lagrangian maps is reduced to 
the transversality of finite dimensional jets of extensions of their reticular Lagrangian embeddings.


%
\section{Weak Caustic-equivalence}
\quad
There exist modalities in the classification list of Section \ref{caust:sec}. 
This means that the caustic-equivalence is still too strong for a generic classification of caustics on a corner.
In order to obtain the generic classification, we need to admit the following equivalence relation:

We say that two function germs in ${\mathfrak M}(r;k+n)^2$ are  {\em weakly reticular ${\cal P}$-${\cal C}$-equivalent}
if they are generating families of  weakly caustic-equivalent reticular Lagrangian maps.
We define the {\em stable}  weakly reticular ${\cal P}$-${\cal C}$-equivalence by the ordinary way. \\

We say that a reticular Lagrangian map $\pi\circ i$ is {\em weakly caustic-stable} if 
the following condition holds:
For any extension $S\in S(T^*{\mathbb R}^n,0)$ of $i$ and any representative 
$\tilde{S}\in S(U,T^*{\mathbb R}^n)$ of $S$, there exists a neighborhood $N_{{\tilde{S}}}$ of 
$\tilde{S}$ such that for any $\tilde{S}'\in N_{{\tilde{S}}}$ the reticular Lagrangian map
$\pi\circ \tilde{S}'|_{\mathbb L}$ at $x_0$ and $\pi\circ i$ are weakly caustic-equivalent for some 
$x_0=(0,\ldots,0,p^0_{r+1},\ldots,p^0_n)\in U$.\\

%
We say that a function germ $F(x,y,u)\in {\mathfrak M}(r;k+n)$ is {\em weakly 
reticular ${\cal P}$-${\cal C}$-stable} if the following condition holds: For any neighborhood $U$ of
$0$ in ${\mathbb R}^{r+k+n}$ and any representative $\tilde{F} \in C^\infty 
(U,{\mathbb R})$ of $F$, 
there exists a neighborhood $N_{\tilde{F}}$ of $\tilde{F}$ in $C^\infty$-topology
 such that for any element $\tilde{G} \in N_{\tilde{F}}$
the germ $\tilde{G}|_{{\mathbb H}^r\times {\mathbb R}^{k+n}}$ at 
$(0,y_0,q_0)$ is weak reticular ${\cal P}$-${\cal C}$-equivalent to $F$
for some $(0,y_0,q_0)\in U$.\\

We here investigate the reticular ${\cal C}$-equivalence classes $B_{2,2,a}^{+,+,2}$ of function germs.
The same methods are valid for the classes $B_{2,2,a}^{\pm,+,1}$, $B_{2,2,a}^{\pm,+,2}$, $B_{2,2,a}^{\pm,-}$.
So we prove only to the classes $B_{2,2,a}^{+,+,2}$.\\

We consider the reticular Lagrangian maps $\pi\circ i_a:
({\mathbb L},0)\rightarrow (T^* {\mathbb R}^2,0) \rightarrow
 ({\mathbb R}^2,0)$ with the generating families
$F_a(x_1,x_2,q_1,q_2)=x_1^2+x_1x_2+ax_2^2+q_1x_1+q_2x_2\ (a>\frac14)$.
We give the caustic of $\pi\circ i_a$ and $\pi\circ i_b$ for $\frac{1}{4}<a<b$.
\begin{figure}[ht]  \begin{center}
\begin{minipage}{0.4\hsize}
  \begin{center}
\includegraphics[width=4.5cm,height=4.5cm]{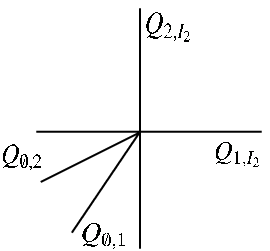} 
\caption{the caustics of $\pi\circ i_a$} 
\end{center}
 \end{minipage}
 \begin{minipage}{0.4\hsize}
  \begin{center}
\includegraphics[width=4.5cm,height=4.5cm]{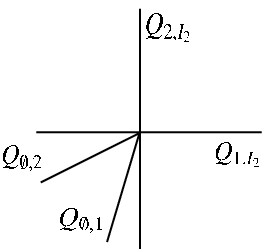} 
\caption{the caustics of $\pi\circ i_b$} 
\end{center}
 \end{minipage}\end{center}
\end{figure}
In these figures $Q_{1,I_2},Q_{2,I_2},Q_{\emptyset,2}$ are in 
the same positions.
Suppose that there exists a diffeomorphism germ $g$ on 
$({\mathbb R}^2,0)$ such that 
$Q_{1,I_2},Q_{2,I_2},Q_{\emptyset,2}$ are invariant under $g$.
Then $g$ can not map $Q_{\emptyset,1}$ from one to the other.
This implies that caustic-equivalence is too strong for generic classifications.
But these caustic are equivalent under the weak caustic-equivalence.
This implies that the reticular Lagrangian map $\pi\circ i_a$ 
is weakly caustic equivalent to $\pi\circ i_1$ for any $a>\frac14$ and hence 
$F_a$ is weakly reticular ${\cal P}$-${\cal C}$-equivalent to $F_1$.
We remark that a homeomorphism germ $g_a$, which gives the weak caustic-equivalence 
of $\pi\circ i_1$ and $\pi\circ i_a$, may be chosen to be 
smooth outside $0$ and depends smoothly on $a$.
This means that the weak caustic-equivalence relation is naturally extended for the
(caustic) stable reticular Lagrangian maps with the generating families
$F'_a(x_1,x_2,q_1,q_2,q_3)=x_1^2+x_1x_2+ax_2^2+q_1x_1+q_2x_2+q_3x_2^2$
and $F'_a$ is weakly reticular ${\cal P}$-${\cal C}$-equivalent to 
$F'(x_1,x_2,q_1,q_2,q_3)=x_1^2+x_1x_2+x_2^2+q_1x_1+q_2x_2$.
The figure of the corresponding caustic is given in \cite[p.602 $B_{2,2}^{+,+,\stackrel<\alpha}$]{retLag}.
We also remark that
the functions 
$x_1^2+x_1x_2+x_2^2+q_1x_1+q_2x_2$ and $x_1^2+x_1x_2+\frac15 x_2^2+q_1x_1+q_2x_2$
are not weakly reticular ${\cal P}$-${\cal C}$-equivalent because $Q_{\emptyset,1}$ and
$Q_{\emptyset,1}$ of  their caustics are in the opposite positions to each other.\\

%
%
%
Therefore  we have that the function germ $f_a(x)=x_1^2+x_1x_2+ax_2^2(a>\frac14)$ are all weak
reticular ${\cal C}$-equivalent. 
Since $\frac{d f_a}{d a}=x^2$ is not included in $\langle x\frac{\partial f_a}{\partial x}\rangle_{{\cal E}(x)}$,
it follows that the $l$-jets of the weak reticular ${\cal C}$-equivalence class of $f_a$ consists
 an immersed manifold 
of $J^l(2,1)$ for $l\geq 2$.

We classify function germs in ${\mathfrak M}(2;k)^2$ with respect to 
the weak reticular ${\cal C}$-equivalence with the codimension$\leq 4$.
Then we have the following list:
\begin{center}
\begin{tabular}{lllll}
\hline $k$ & Normal form & codim & Notation \\ 
\hline
 $0$ & $x_1^2\pm x_1x_2+ \frac15 x_2^2$ & $2$  &  $B_{2,2}^{\pm,+,1}$\\
 & $x_1^2\pm x_1x_2+ x_2^2$ & $2$ &   $B_{2,2}^{\pm,+,2}$\\
 & $x_1^2\pm x_1x_2-x_2^2$ & $2$ & $B_{2,2}^{\pm,-}$\\
 & $x_1^2\pm x_2^2$ & $3$ &      $B_{2,2}^{\pm,0}$\\
 & $(x_1\pm x_2)^2\pm x_2^3$ & $3$ &  $B_{2,2,3}^{\pm,\pm}$\\
& $x_1^2\pm x_1x_2\pm x_2^3$ & $3$  &
                              $B_{2,3}^{\pm,\pm}$\\
& $x_1^2\pm x_1x_2^2\pm  x_2^3$ & $4$ & 
                              $B_{2,3'}^{\pm,\pm}$ \\
 & $x_1^3\pm x_1^2x_2\pm  x_2^2$ & $4$ & 
                              $B_{3,2'}^{\pm,\pm}$ \\ 
 & $x_1^3\pm x_1x_2\pm  x_2^2$ & $3$ &
                              $B_{3,2}^{\pm,\pm}$ \\
\hline 
 $1$ & $\pm y^3+ x_1y\pm  x_2y+x_2^2$ & $3$ 
       & $C_{3,2}^{\pm,\pm}$ \\
       & $\pm y^3+ x_1y\pm  x_2y^2+x_2^2$ & $4$ 
  & $C_{3,2,1}^{\pm,\pm}$ \\
  & $\pm y^3+ x_2y\pm  x_1y^2+x_1^2$ & $4$ 
      
   & $C_{3,2,2}^{\pm,\pm}$ \\
\hline 
\end{tabular}
\end{center}

\begin{prop}
Let $\pi\circ i_a:({\mathbb L},0)\rightarrow (T^* {\mathbb R}^2,0) \rightarrow
 ({\mathbb R}^2,0)$ be the reticular Lagrangian map with the generating family
$x_1^2+x_1x_2+ax_2^2+q_1x_1+q_2x_2$.
Let $S_a\in S(T^* {\mathbb R}^2,0)$ be extensions of $i_a$.
Then the weak caustic-equivalence class 
\[ [j^lS_1(0)]_w:=\bigcup_{a>\frac{1}{4}}[j^lS_a(0)]_c \]
is an immersed manifold in $S^l(2)$ for $l\geq 1$.
\end{prop}
{\em Proof}. 
The function germ $G_a(Q_1,Q_2,q_1,q_2)=Q_1^2+Q_1Q_2+aQ_2^2+q_1Q_1+q_2Q_2$
is a generating function of the canonical relation $P_{S_a}$ and we have that 
\[ S_a(Q,P)=(-(2Q_1+Q_2+P_1),-(Q_1+2aQ_2+P_2),Q_1,Q_2).\]
This means that $\frac{d S_a}{da}=(0,-2Q_2,0,0)=X_f\circ S_a$ for
 $f=-p_2^2$.
Suppose that $j^1(\frac{d S_a}{da})(0)$ is included in $T_z(rLa^1(2)\cdot z)$.
Then there exist  $h_1,h_2\in {\mathfrak M}_{Q,P}$ and $h_0\in {\mathfrak M}_{Q,P}^2$ such that
\begin{eqnarray*}
-Q_2^2 \equiv h_1(q\circ S_a)Q_1+
h_2(q\circ S_a)Q_2+
h_0(q\circ S_a) \mbox{ mod }  \langle Q_1P_1,Q_2P_2 \rangle_{{\cal E}_{Q,P}}+
{\mathfrak M}_{Q,P}^3.
\end{eqnarray*}
We need only to consider the linear parts of $h_1,h_2$ and the quadratic part of $h_0$.
The coefficients of $P_1^2,P_2^2,P_1P_2$ depend only on
the coefficients of $Q_1^2,Q_2^2,Q_1Q_2$ on $h_0$ respectively.
This means that $h_0\equiv 0$.
We set $h_1=bq_1+cq_2,h_2=dq_1+eq_2$ and 
calculate the coefficients of $Q_1^2,Q_1Q_2,Q_1P_2,Q_2P_1$ in the equation.
Then we have that $-2 b - c = 0,
-b - 2 d - e - 2 c a =0,c=0,d=0$.
This means that $e=0$.
Then we have that the coefficient of $Q_2^2$ of the right hand side of the equation
is equivalent to $-d-ae=0$.
This contradicts the equation. \hfill $\blacksquare$\\

We consider the (caustic) stable reticular Lagrangian map
$\pi\circ i_a:({\mathbb L},0)\rightarrow (T^* {\mathbb R}^3,0) \rightarrow
 ({\mathbb R}^3,0)$ with the generating family
 $x_1^2+x_1x_2+ax_2^2+q_1x_1+q_2x_2+q_3 x_2^2$ and 
 take an extension $S'_a\in S(T^* {\mathbb R}^2,0)$ of $i_a$,
 then we have by the analogous method that:
\begin{coro}
Let $S'_a$ be as above.
Then the weak caustic-equivalence class
\[ [j^lS'_1(0)]_w:=\bigcup_{a>\frac{1}{4}}[j^lS'_a(0)]_c \]
is an immersed manifold in $S^l(3)$ for $l\geq 1$.
\end{coro}

Since the caustic of $\pi\circ i_a$ is given by the restrictions of $\pi\circ i_a$
to $L^0_\sigma\cap L^0_\tau$ for $\sigma\neq \tau$ in this case,
it follows that the caustic is determined by the linear part of $i_a$.
This means that $\pi\circ i_a$ is $1$-determined with respect to  the weak caustic-equivalence
(cf., Definition \ref{caustequi:dfn}). 
\begin{tth}
The function germ $F(x_1,x_2,q_1,q_2)=x_1^2+x_1x_2+x_2^2+q_1x_1+q_2x_2$ is 
a weakly reticular ${\cal P}$-${\cal C}$-stable unfolding of $f(x_1,x_2)=x_1^2+x_1x_2+x_2^2$
\end{tth}
{\em Proof}.
We define $F'\in {\mathfrak M}(2;3)^2$ by $F'(x_1,x_2,q_1,q_2,q_3)=F(x_1,x_2,q_1,q_2)+q_3 x_2^2$
Then $F'$ is a reticular ${\cal P}$-${\cal R}^+$-stable unfolding of $f$.
It follows that for any neighborhood $U'$ of $0$ in ${\mathbb R}^5$ and
any representative $\tilde{F'}\in C^\infty(U,{\mathbb R})$,
there exists a neighborhood $N_{\tilde{F'}}$ such that 
for any $\tilde{G'}\in N_{\tilde{F'}}$ 
the function germ $\tilde{G'}|_{{\mathbb H}^2\times {\mathbb R}^3}$ at $p'_0$ is reticular ${\cal P}$-${\cal R}^+$-equivalent to 
$F'$ for some $p'_0=(0,0,q^0_1,q^0_2,q^0_3)\in U'$.

Let a neighborhood $U$ of $0$ in ${\mathbb R}^4$ and a
representative $\tilde{F}\in C^\infty(U,{\mathbb R})$ be given.
We set the open interval $I=(-0.5,0.5)$ and set $U'=U\times I$.
Then there exists $N_{\tilde{F'}}$ for which the above condition holds.
We can choose a neighborhood $N_{\tilde{F}}$ of $\tilde{F}$ such that
for any $\tilde{G}\in N_{\tilde{F}}$ the function $\tilde{G}+q_3x_2^2\in N_{\tilde{F'}}$.
Let a function $\tilde{G}\in N_{\tilde{F}}$ be given.
Then the function germ $G'=(\tilde{G}+q_3x_2^2)|_{{\mathbb H}^2\times {\mathbb R}^3}$ at $p'_0$
is reticular ${\cal P}$-${\cal R}^+$-equivalent to 
$F'$ for some $p'_0=(0,0,q_1^0,q_2^0,q_3^0)\in U'$.
We define $G\in {\mathfrak M}(2;2)^2$ by $\tilde{G}|_{{\mathbb H}^2\times {\mathbb R}^2}$ at 
$p_0=(0,0,q_1^0,q_2^0)\in U$.
Then it holds that $G'(x,q)=G(x,q_1,q_2)+(q_3+q_3^0)x_2^2$, and
$G'|_{q=0}=G(x,0)+q_3^0x_2^2$ is reticular ${\cal R}$-equivalent to $f$.
Let $(\Phi,a)$ be the reticular ${\cal P}$-${\cal R}^+$-equivalence from $G'$ to $F'$.
We write $\Phi(x,q)=(x\phi_1(x,q),\phi_1^2(q),\phi_2^2(q),\phi_3^2(q))$.
By shrinking $U$ if necessary, we may assume that the map germ
\[ (q_1,q_2)\mapsto (\phi_1^2(q_1,q_2,0),\phi_2^2(q_1,q_2,0)) \mbox{ on } ({\mathbb R}^2,0) \]
is a diffeomorphism germ.
Then $F$ is  reticular ${\cal P}$-${\cal R}^+$-equivalent to
$G_1\in {\mathfrak M}(2;2)^2$ given by $G_1(x,q)=G(x_1,x_2,q_1,q_2)+(\phi^2_3(q_1,q_2,0)+q^0_3)x_2^2$.
It follows that the reticular Lagrangian maps defined by $F$ and $G_1$ are Lagrangian equivalent.
We have that 
\[  j^2(G+q_3^0x_2^2)(0)=j^2G_1(0),\ \ q_3^0>-0.5.\]
This means that the caustic of $G_1$ is weakly caustic-equivalent to the caustic of $G$ 
because the reticular Lagrangian maps of  $G_1$ and $F$ are the same weak caustic-equivalence class  that is
$1$-determined under the weak caustic-equivalence.
This means that $F$ and $G$ are weakly reticular ${\cal P}$-${\cal C}$-equivalent.
Therefore $F$ is weakly reticular ${\cal P}$-${\cal C}$-stable.
\hfill $\blacksquare$\\

By the above consideration, we have that: 
For each singularity $B_{2,2}^{\pm,+,1}, B_{2,2}^{\pm,+,2}, B_{2,2}^{\pm,-}$,
if we take the symplectic diffeomorphism germ $S_a(S'_a)$ as the above method, then
the weak caustic-equivalence class $[j^lS_a(0)]_w([j^lS'_a(0)]_w)$ is one class and  immersed manifold in 
$S^l(2)(S^l(3))$ for 
$l\geq 1$ respectively.\\

We now start to prove the main theorem: 
We choose the weakly caustic-stable reticular Lagrangian maps $\pi\circ i_X:
({\mathbb L},0)\rightarrow (T^* {\mathbb R}^n,0)
\rightarrow  ({\mathbb R}^n,0)$ 
for 
\begin{equation}
 X=B_{2,2}^{\pm,+,1},B_{2,2}^{\pm,+,2},
B_{2,2}^{\pm,-}. \label{weak:X}
\end{equation}
We also choose the caustic-stable reticular Lagrangian maps $\pi\circ i_X:
({\mathbb L},0)\rightarrow (T^* {\mathbb R}^3,0)
\rightarrow  ({\mathbb R}^3,0)$ for 
\begin{equation}
X=B_{2,2}^{\pm,0},B_{2,2,3}^{\pm.\pm},B_{2,3}^{\pm,\pm},
B_{3,2}^{\pm,\pm},B_{2,3'}^{\pm,\pm},B_{3,2'}^{\pm,\pm},
C_{2,3}^{\pm,\pm},C_{3,2,1}^{\pm,\pm},C_{3,2,2}^{\pm,\pm}.\label{caust:X}
\end{equation}
Then other reticular Lagrangian maps are not caustic-stable since other singularities have 
reticular ${\cal C}$-codimension $>4$.
We choose extensions $S_X\in S(T^* {\mathbb R}^n,0)$ 
of $i_X$ for all $X$.
We define that
\[ O'_1=\{ \tilde{S}\in S(U,T^* {\mathbb R}^n)\ |\ 
j^{n+2}_0\tilde{S} \mbox{ is transversal to }[j^{n+2}S_X(0)]_w
\mbox{ for all } X \mbox{ in } (\ref{weak:X}) \},\]
\[ O'_2=\{ \tilde{S}\in S(U,T^* {\mathbb R}^n)\ |\ 
j^{n+2}_0\tilde{S} \mbox{ is transversal to }[j^{n+2}S_X(0)]_c
\mbox{ for all } X \mbox{ in } (\ref{caust:X}) \},\]
where $j^l_0\tilde{S}(x)=j^l\tilde{S}_x(0)$.
Then $O'_1$ and $O'_2$ are residual sets. We set 
\[ Y=\{ j^{n+2}S(0)\in S^{n+2}(n)\ |\ 
\mbox{the codimension of }[j^{n+2}S(0)]_L>10\}.\]
Then $Y$ is an algebraic set in $S^{n+2}(n)$ by \cite[Theorem 6.6 (a')]{generic}.
Therefore we can define that 
\[ O''=\{ \tilde{S}\in S(U,T^* {\mathbb R}^n)\ |\ 
j^{n+2}_0\tilde{S} \mbox{ is transversal to } Y \}.\]  
For any $S\in S(T^* {\mathbb R}^n,0)$ with $j^{n+2}S(0)$ and any generating family $F$ of $\pi\circ S|_{{\mathbb L}}$,
the function germ $F|_{q=0}$ has the reticular ${\cal R}^+$-codimension $>5$.
This means that $F|_{q=0}$ has the reticular $C$-codimension $>4$.
It follows that $j^{n+2}S(0)$ does not belong to the above equivalence classes.
Then $Y$ has codimension $>8$.
Then we have that
\[ O''=\{ \tilde{S}\in S(U,T^* {\mathbb R}^n)\ |\ 
j^{n+2}_0\tilde{S}(U)\cap Y=\emptyset \}.\]
We define $O=O'_1\cap O'_2\cap O''$. 
Since all $\pi\circ i_X$ for $X$ in (\ref{weak:X}) are weak caustic $1$-determined, and
all $\pi\circ i_X$ in (\ref{caust:X}) are caustic $5$-determined by Lemma \ref{findet:lem}.
Then $O$ has the required condition.
\hfill $\blacksquare$
\begin{figure}[ht]\begin{center}
$B_{2,2}^{*,*}$ in $2D$\\
 \begin{minipage}{0.30\hsize}
  \begin{center}
    \includegraphics[width=4cm,height=4cm]{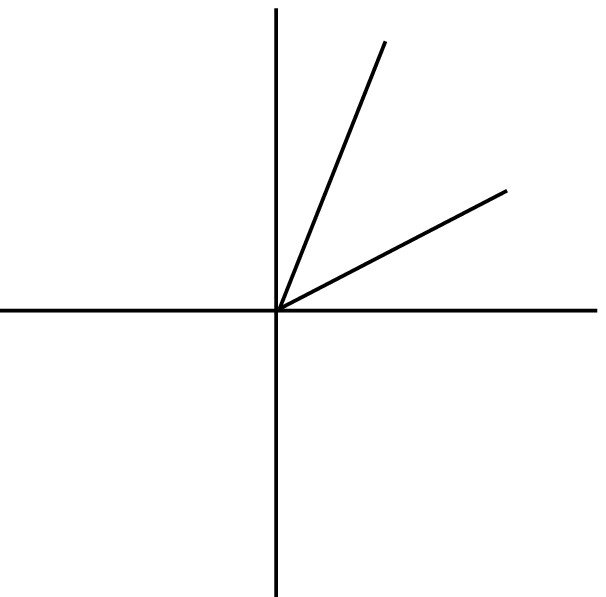}
  \caption{$B_{2,2}^{+,+,1},\ B_{2,2}^{+,+,2}$} \end{center}
 \end{minipage}
\begin{minipage}{0.30\hsize}
  \begin{center}
    \includegraphics[width=4cm,height=4cm]{B22-+.eps} 
  \caption{$B_{2,2}^{-,+,1},\ B_{2,2}^{-,+,2}$}\end{center}
 \end{minipage}
 \begin{minipage}{0.30\hsize}
  \begin{center}
    \includegraphics[width=4cm,height=4cm]{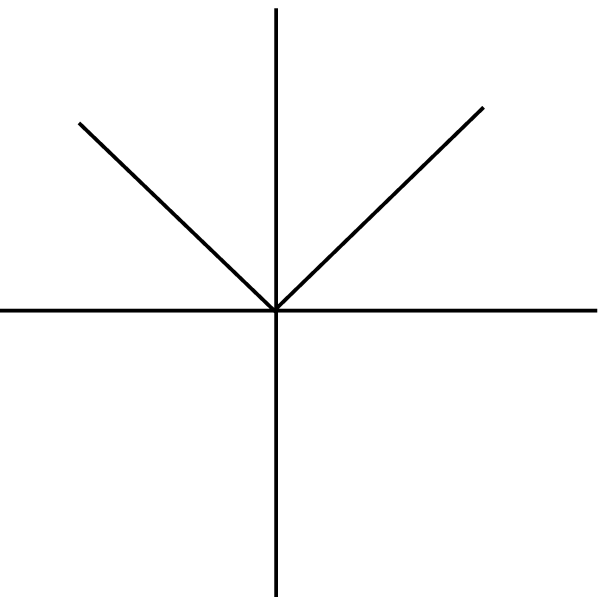} 
  \caption{$B_{2,2}^{+,-},\ B_{2,2}^{-,-}$} \end{center}
 \end{minipage}\end{center}
\end{figure}
\begin{figure}[ht]\begin{center}
$B_{2,2}^{*,*}$ in $3D$\\
 \begin{minipage}{0.30\hsize}
  \begin{center}
    \includegraphics[scale=1,clip]{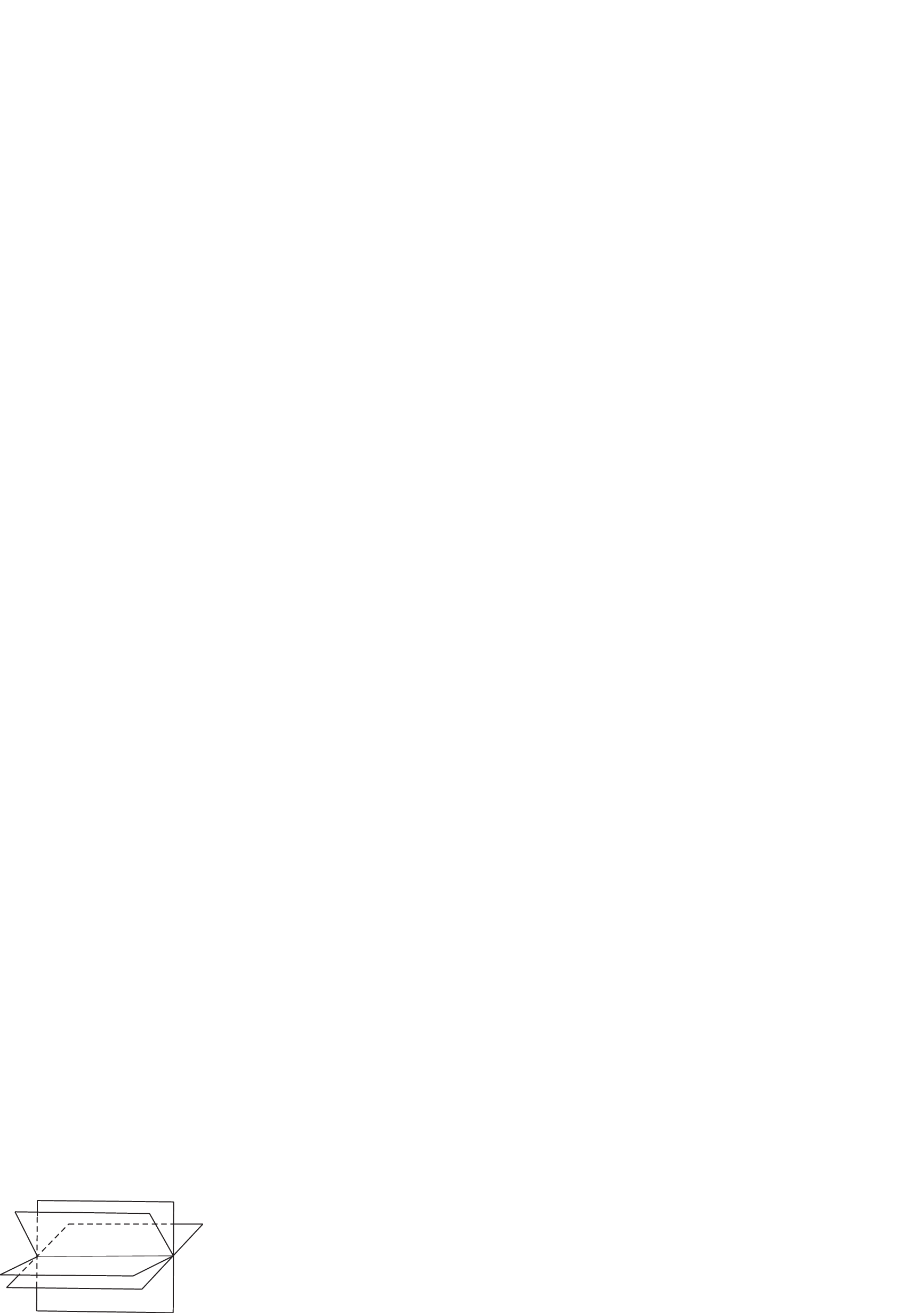}
  \caption{$B_{2,2}^{+,+,1},\ B_{2,2}^{+,+,2}$} \end{center}
 \end{minipage}
\begin{minipage}{0.30\hsize}
  \begin{center}
    \includegraphics[scale=1]{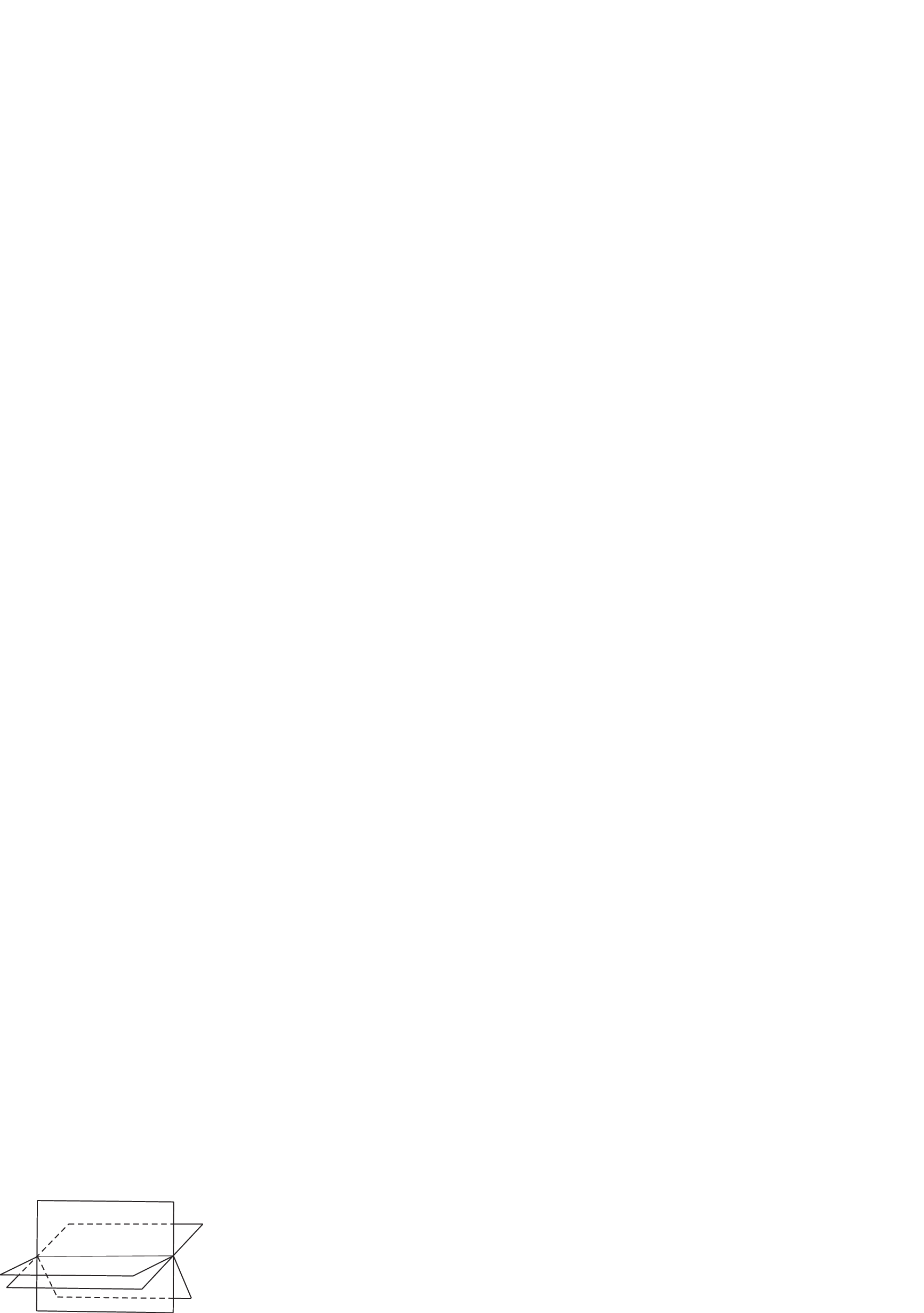} 
  \caption{$B_{2,2}^{-,+,1},\ B_{2,2}^{-,+,2}$}\end{center}
 \end{minipage} 
 \begin{minipage}{0.30\hsize}
  \begin{center}
    \includegraphics[scale=1]{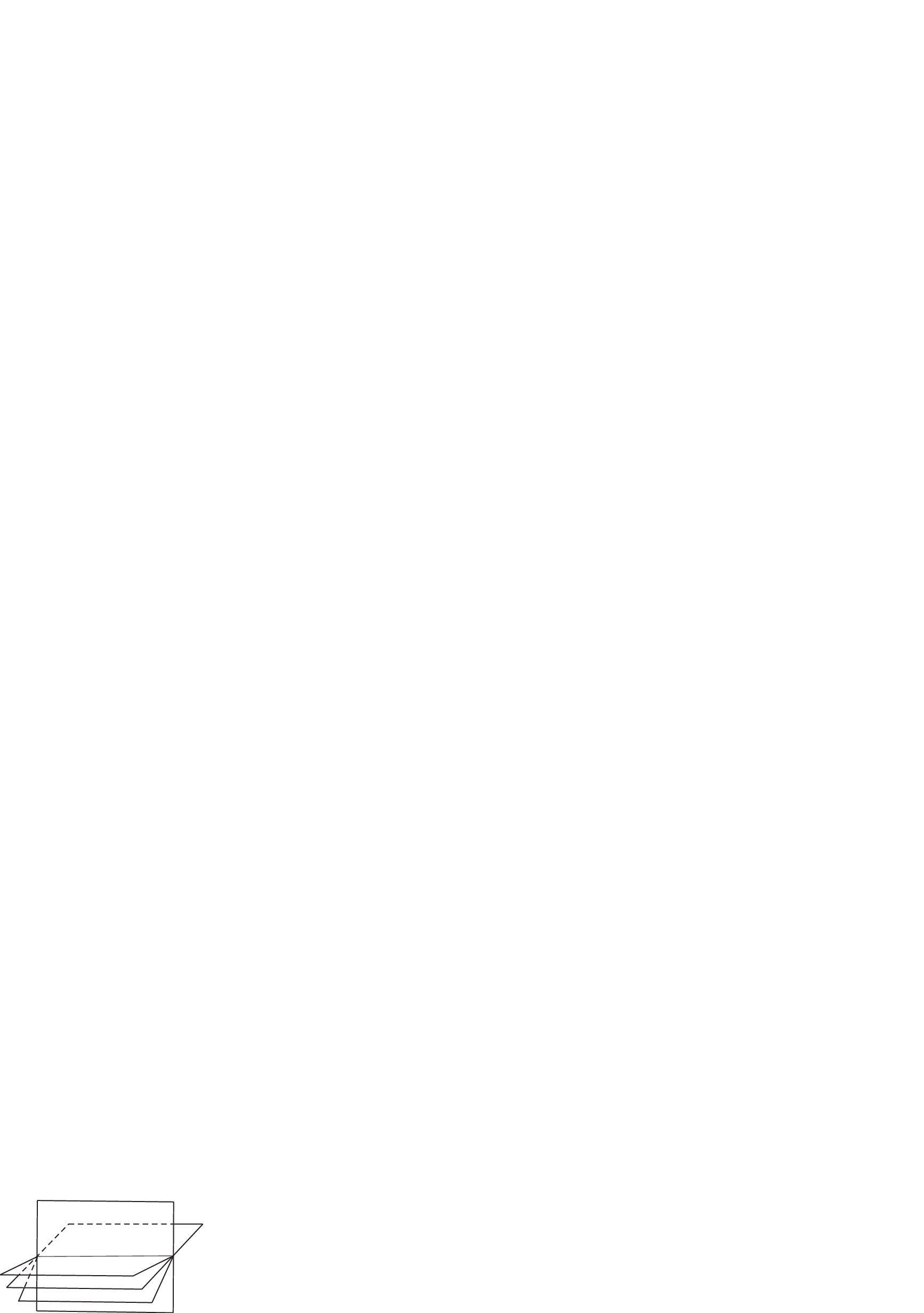} 
  \caption{$B_{2,2}^{+,-},\ B_{2,2}^{-,-}$} \end{center}
 \end{minipage}\end{center}
\end{figure}
\begin{figure}[ht]\begin{center}
 \begin{minipage}{0.30\hsize}
  \begin{center}
    \includegraphics[width=5.5cm,height=5cm]{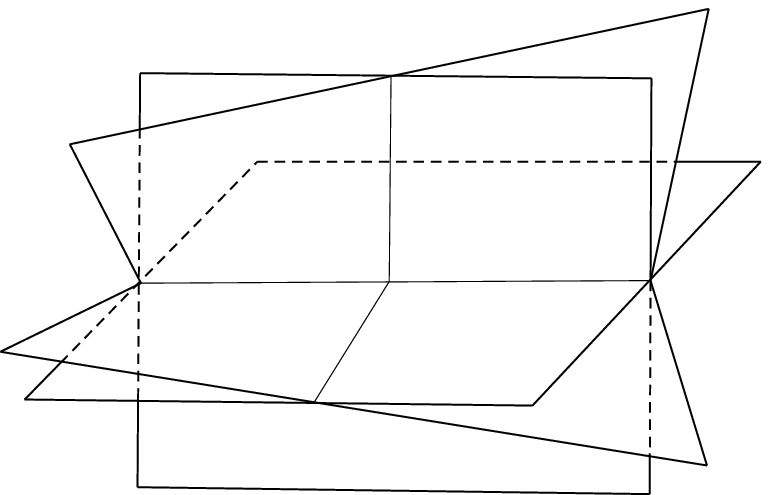}
  \caption{$B_{2,2}^{+,0}$}  \end{center}
 \end{minipage}\hspace{2cm}
 \begin{minipage}{0.30\hsize}
  \begin{center}
    \includegraphics[width=5.5cm,height=5cm]{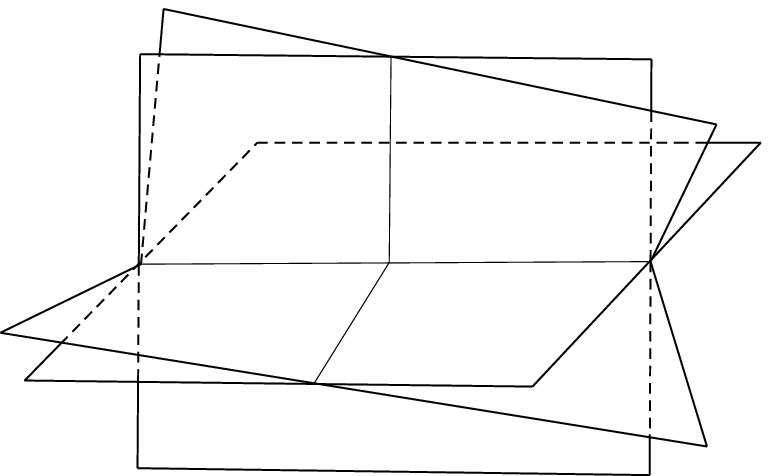}
  \caption{$B_{2,2}^{-,0}$}\end{center}
 \end{minipage}\end{center}
\end{figure}
\begin{figure}[ht]\begin{center}
 \begin{minipage}{0.30\hsize}
  \begin{center}
    \includegraphics[width=5cm,height=5cm]{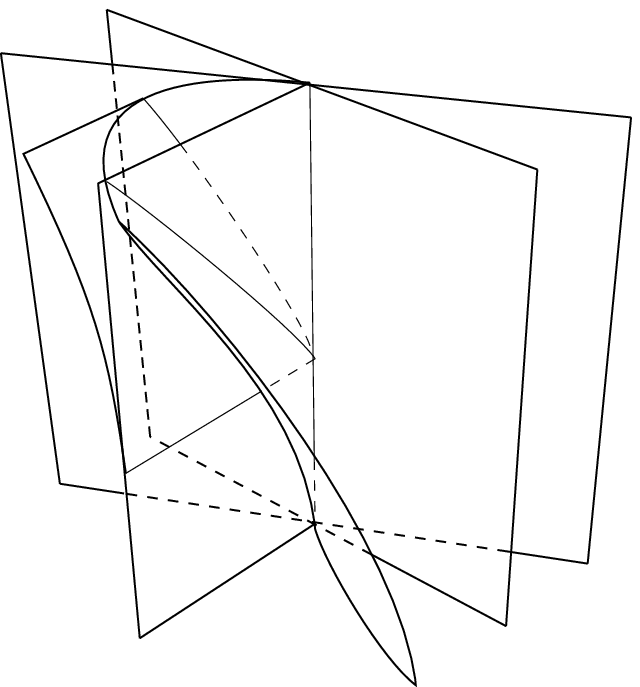} 
  \caption{$B_{2,2,3}^{+,+}$}\end{center}
 \end{minipage}\hspace{2cm}
 \begin{minipage}{0.30\hsize}
  \begin{center}
    \includegraphics[width=5cm,height=5cm]{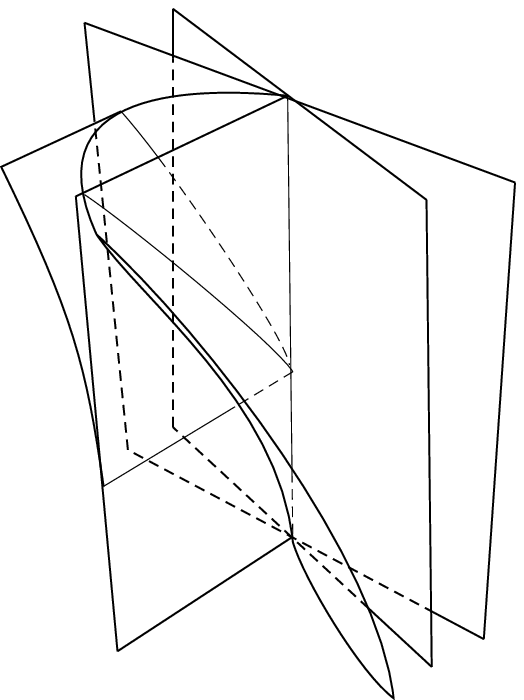}
  \caption{$B_{2,2,3}^{+,-}$}  \end{center}
 \end{minipage}\vspace{5mm}\\
  \begin{minipage}{0.30\hsize}
  \begin{center}
    \includegraphics[width=4.5cm,height=5cm]{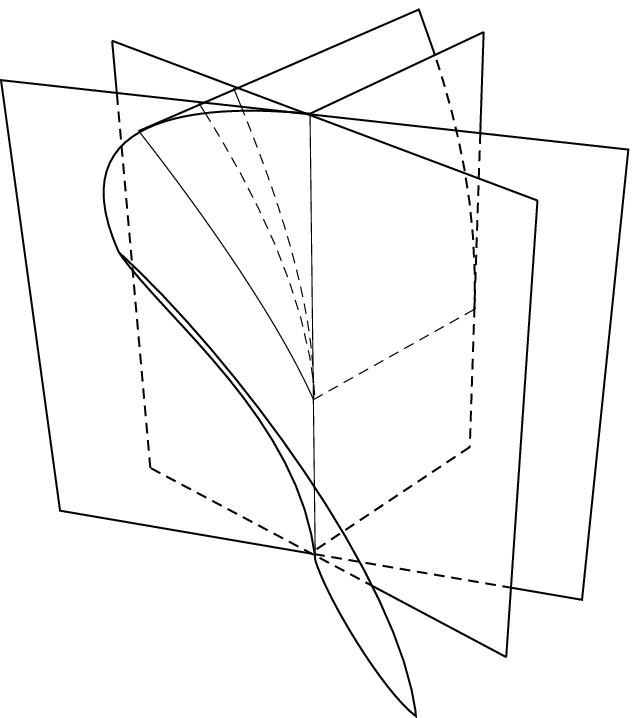}
  \caption{$B_{2,2,3}^{-,+}$}  \end{center}
 \end{minipage}\hspace{2cm}
 \begin{minipage}{0.30\hsize}
  \begin{center}
    \includegraphics[width=4cm,height=5cm]{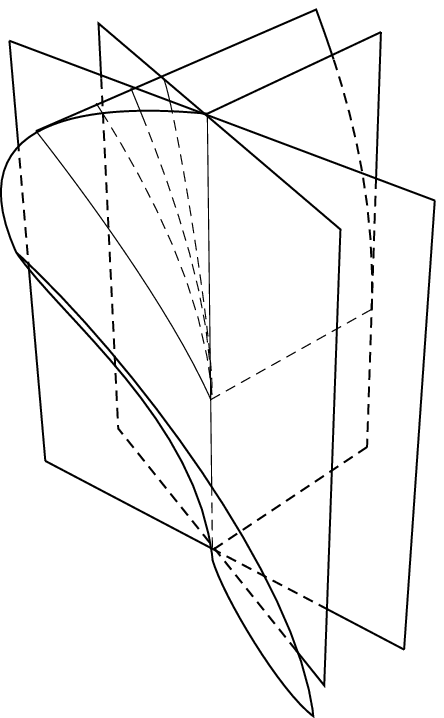}  
  \caption{$B_{2,2,3}^{-,-}$}\end{center}
 \end{minipage}\end{center}
\end{figure}
\begin{figure}[ht]\begin{center}
 \begin{minipage}{0.30\hsize}
  \begin{center}
    \includegraphics[width=5cm,height=5cm]{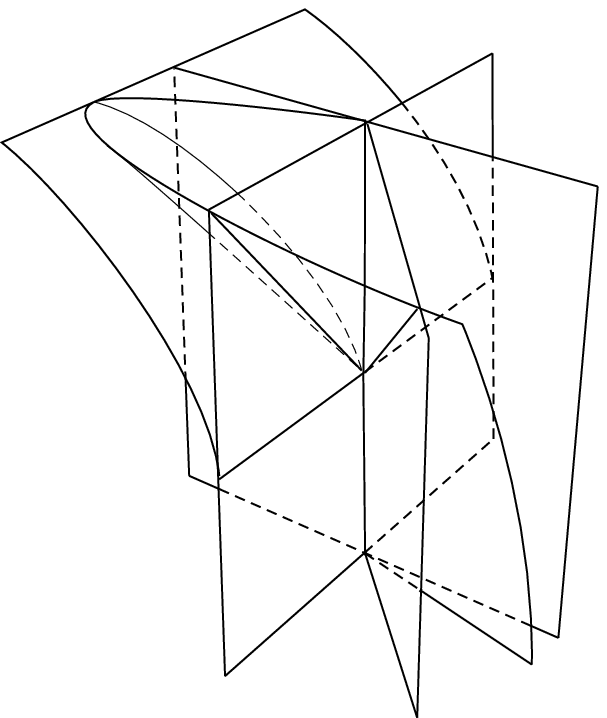}
  \caption{$B_{2,3}^{+,+}$}
    \end{center}
 \end{minipage}\hspace{2cm}
 \begin{minipage}{0.30\hsize}
  \begin{center}
    \includegraphics[width=5cm,height=5cm]{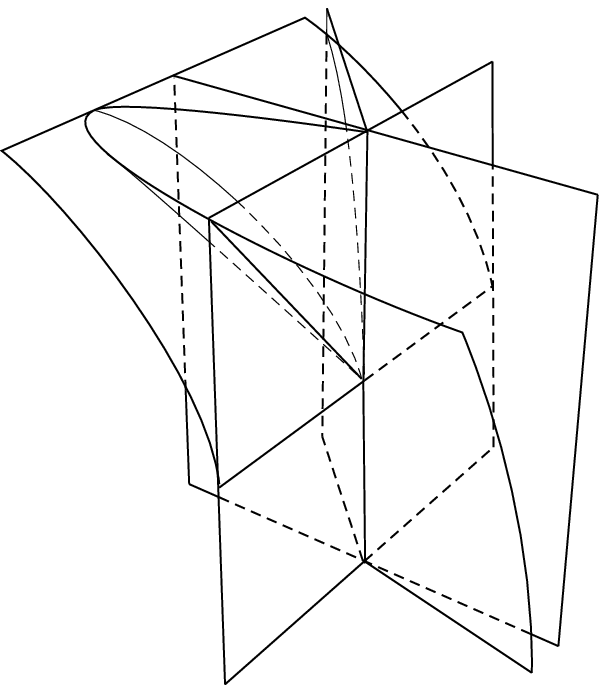}
  \caption{$B_{2,3}^{+,-}$}
  \end{center}
 \end{minipage}\vspace{5mm}\\
\begin{minipage}{0.30\hsize}
  \begin{center}
    \includegraphics[width=5cm,height=5cm]{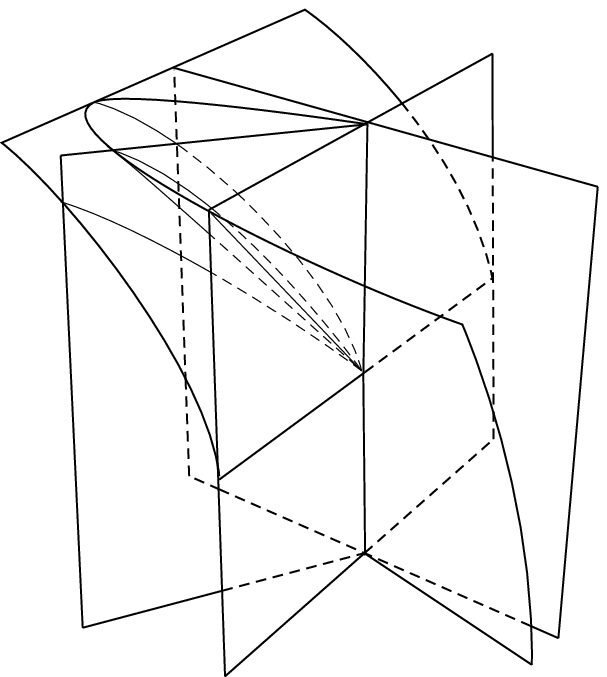}  
  \caption{$B_{2,3}^{-,+}$}\end{center}
 \end{minipage}\hspace{2cm}
 \begin{minipage}{0.30\hsize}
  \begin{center}
    \includegraphics[width=5cm,height=5cm]{B23--.eps}
  \caption{$B_{2,3}^{-,-}$}  \end{center}
 \end{minipage}\end{center}
\end{figure}
\begin{figure}[ht]\begin{center}
 \begin{minipage}{0.30\hsize}
  \begin{center}
    \includegraphics[width=5cm,height=5cm]{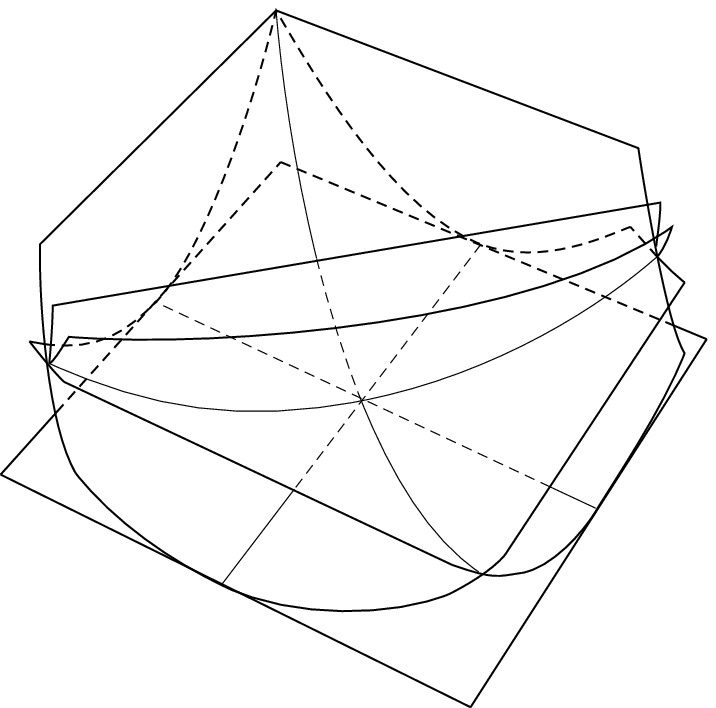}
  \caption{$C_{3,2}^{+,+}$} \end{center}
 \end{minipage}\hspace{2cm}
 \begin{minipage}{0.30\hsize}
  \begin{center}
    \includegraphics[width=5cm,height=5cm]{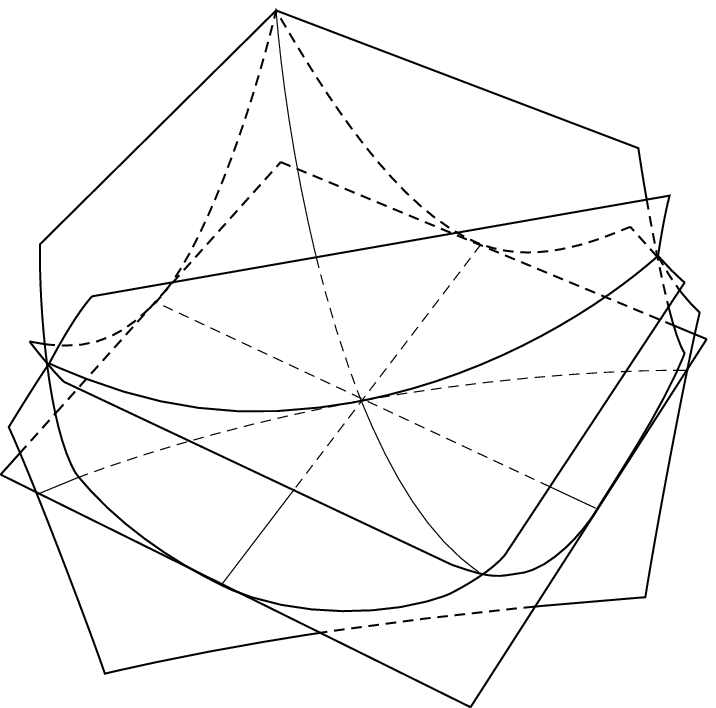}
  \caption{$C_{3,2}^{+,-}$}\end{center}
 \end{minipage}\vspace{5mm}\\
\begin{minipage}{0.30\hsize}
  \begin{center}
    \includegraphics[width=5cm,height=5cm]{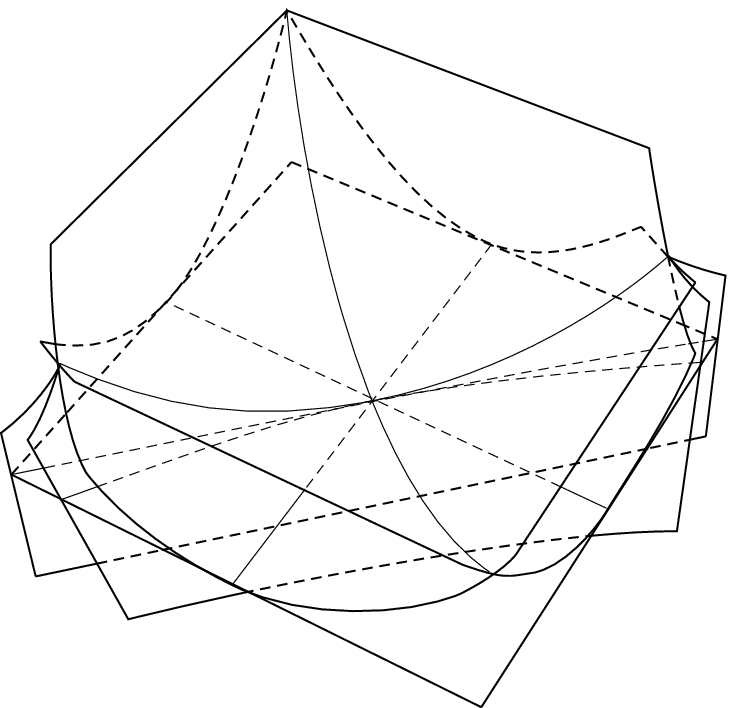}
  \caption{$C_{3,2}^{-,+}$}\end{center}
 \end{minipage}\hspace{2cm}
 \begin{minipage}{0.30\hsize}
  \begin{center}
    \includegraphics[width=5cm,height=5cm]{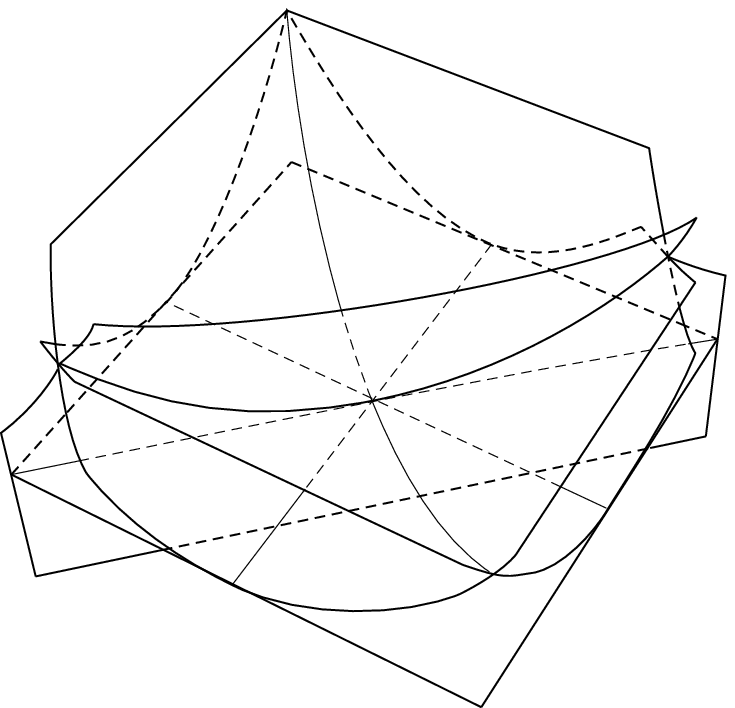}
  \caption{$C_{3,2}^{-,-}$}\end{center}
 \end{minipage}\end{center}
\end{figure}

\bibliographystyle{plain}
\bibliography{bibt}
\end{document}